\documentclass[A4,twoside,11pt,leqno]{article}
\usepackage{amssymb}
\usepackage{amsthm}
\usepackage[tbtags]{amsmath}
\usepackage{doc}
\usepackage{latexsym}
\usepackage{amscd}
\usepackage[frame,cmtip,arrow,matrix,line,graph,curve]{xy}
\usepackage{graphics}
\usepackage{epsfig}
\usepackage{eucal}
\usepackage{mathrsfs}
\font\headd=cmr8
\newtheorem{theorem}{Theorem}
\newtheorem*{thm}{Theorem}
\newtheorem{lemma}{Lemma}

\newtheorem{conjecture}{Conjecture}
\newtheorem{corollary}{Corollary}

\def\Z {\mathbb{Z}}
\def\Q {\mathbb{Q}}
\def\R {\mathbb{R}}
\def\C {\mathbb{C}}
\def\QQ {\overline{\Q}}

\def\u {\mathcal{U}}

\begin{document}
\thispagestyle{plain}
 \markboth{}{}
\small{\addtocounter{page}{0} \pagestyle{plain}
\noindent{\scriptsize KYUNGPOOK Math. J. 00(0000), 000-000}
\vspace{0.2in}\\
\noindent{\large\bf Some transcendental functions with an empty exceptional set}
\footnote{{}\\ \\[-0.7cm]
* Corresponding Author.\\
2000 Mathematics Subject Classification: Primary 11J81 .\\
Key words and phrases: Hermite-Lindemann theorem, Lindemann-Weierstrass theorem, Baker theorem, Mahler classification, Schanuel conjecture .\\
}
\vspace{0.15in}\\

\noindent{\sc F. M. S. Lima}
\newline
{\it Institute of Physics, University of Brasilia, Brasilia, DF, Brazil\\
e-mail} : {\verb|fabio@fis.unb.br|}
\vspace{0.15in}\\
\noindent{\sc Diego Marques$^*$}
\newline
{\it Department of Mathematics, University de Brasilia, Brasilia, DF, Brazil\\
e-mail} : {\verb|diego@mat.unb.br|}
\vspace{0.15in}\\
{\footnotesize {\sc Abstract.}
A transcendental function usually returns transcendental values at algebraic points. The (algebraic) exceptions form the so-called \emph{exceptional set}, as for instance the unitary set $\{0\}$ for the function $f(z) = e^z \,$, according to the Hermite-Lindemann theorem. In this note, we give some explicit examples of transcendental entire functions whose exceptional set are empty.
}
\vspace{0.2in}\\
\pagestyle{myheadings}
 \markboth{\headd F. M. S. Lima, D. Marques$~~~~~~~~~~~~~~~~~~~~~~~~~~~~~~~~~~~~~~~~~~~~~\,$}
 {\headd $~~~~~~~~~~~~~~~~~~~~~~~~~~~~~~~~~~~~~$Some transcendental functions with an empty exceptional set}
%
%
%
%

\section{Introduction}

\indent
An algebraic function is a function $f(x)$ which satisfies $P(x,f(x))=0$ for some nonzero polynomial $P(x,y)$ with complex coefficients. Functions that can be constructed using only a finite number of elementary operations are examples of algebraic functions. A function which is not algebraic is, by definition, a \emph{transcendental function} --- e.g., basic trigonometric functions, exponential function, their inverses, etc. If $f$ is an \textit{entire} function, namely a function which is analytic in $\mathbb{C}$, to say that $f$ is a transcendental function amounts to say that it is not a polynomial. By evaluating a transcendental function at an algebraic point of its domain, one usually finds a transcendental number, but exceptions can take place.  For a given transcendental function, the set of all exceptions (i.e., all algebraic numbers of the function domain whose image is an algebraic value) form the so-called \emph{exceptional set} (denoted by $S_f$). This set plays an important role in transcendental number theory (see, e.g., Ref.~\cite{Wald} and references therein). For instance, it can be used for proving that $\,e\,$ and $\,\pi\,$ are transcendental numbers~\cite{Gordan}. The Hermite-Lindemann theorem (1884) was the first general result in this direction~\cite{Baker}. For our purposes, it is suitable to state this theorem in the following simple form.

\begin{thm}[Hermite-Lindemann]
\label{HLtheorem}
\; The number $\,e^z\,$ is transcendental for all non-zero algebraic values of $\,z$.
\end{thm}

Therefore, $S_f=\{0\}$ for $f(z)=e^z$. This `almost empty' result led people to search for a transcendental function that could yield transcendental values for all algebraic points. The existence of such transcendental functions was first proved by St\"{a}ckel (1895), as a corollary of his main theorem in Ref.~\cite{Sta}.

\begin{thm}[St\"{a}ckel]
\; For any countable set $\,A \subseteq \C\,$ and each dense set $\,T \subseteq \C$, there is a transcendental entire function $f$ such that $f(A) \subseteq T$.
\end{thm}

By putting $\,A=\QQ\,$ and $\,T = \C \, \backslash \QQ\,$, where $\QQ\,$ is the set of all complex algebraic numbers, it follows that there exists a transcendental function whose exceptional set is empty.\footnote{There are, of course, several simple examples of \emph{algebraic} functions with empty exceptional sets, e.g. $f(z) = \pi +z$, but we are interested here in \emph{transcendental} functions only.} In 2006, Surroca \cite{Surroca} remarked that under the hypothesis of the Schanuel's conjecture the function $e^{e^z}$ takes transcendental values for all $z \in \QQ$. 

In a recent paper, Huang et al \cite{D3} proved that any $A\subseteq \QQ$ is the exceptional set of uncountable many transcendental entire functions $f$ (see \cite{D4} for the hypertranscendental version).

The aim of this paper is to explicit the constructions made in \cite{D3} in order to produce examples of such functions in the case $A=\emptyset$. Along the way, we find some simple examples of transcendental functions with an empty exceptional set\footnote{We remark that the choice of $A=\emptyset$ is arbitrary and so we encourage the reader to follow our approach with other choices of $A\subseteq \QQ$}. More precisely, our main result is the following


Let $\{\alpha_1,\alpha_2,\ldots\}$ be an enumeration of $\QQ$. We define the function $\u: \C \times \C \rightarrow \C$ as
\begin{equation}
\label{d}
\u(w,z) := \sum_{n=1}^\infty{\frac{w^n}{\left[1+\sum_{k=1}^n|\sigma_k(\alpha_1,\ldots,\alpha_n)|\right] \left(|w|^n+1\right) \, n!}} \left(z-\alpha_1\right) \cdots \left(z-\alpha_n \right),
\end{equation}
where $\sigma_k$ is the $k$-th elementary symmetric polynomial.

\begin{theorem}\label{t2}
For any positive real number $t$ and any algebraic number $\alpha$, the number $\u(t,\alpha)$ is transcendental if and only if $t$ is transcendental. In particular, we have
\begin{center}
$S_{\u(t,\cdot)} = \left\{
\begin{array}{rcl}
\emptyset\,,& \mbox{if} & t\ \mbox{is transcendental}\\
\QQ\,, & \mbox{if} & t\ \mbox{is algebraic} \, .
\end{array}
\right.$
\end{center}
\end{theorem}

\section{Some simple examples}

Let $\,f\!: \C \rightarrow \C\,$ be a transcendental function and let us denote by $S_f$ the \emph{exceptional set} of $f$, i.e. the set of all $z \in \QQ$ for which $f(z) \in \QQ$.  Let us start our search for transcendental entire functions with empty exceptional sets by showing the truth of the Surroca's remark, mentioned in the previous section. For that, we recall the Schanuel's conjecture, one of the main open problems in transcendental number theory.

\begin{conjecture}[Schanuel]
\; If $ z_1, \ldots , z_n$ are complex numbers linearly independent over $\,\Q$, then among the numbers $z_1 , \ldots , z_n , e^{z_1} , \ldots , e^{z_n}$, at least $\,n\,$ are algebraically independent.
\end{conjecture}

This conjecture was introduced in the 1960's by Schanuel in a course given by Lang~\cite{Lang}. It has several important consequences, as for instance: if $\alpha$ is a non-zero algebraic number, the numbers $\alpha,\ e^{\alpha}$ are $\Q$-linearly independent (as consequence of Hermite-Lindemann theorem), then at least two distinct numbers among $\alpha,e^{\alpha},e^{\alpha},e^{e^{\alpha}}$ are transcendental. Therefore $e^{e^{\alpha}}$ is transcendental. The Surroca's result follows by noting that $e^{e^0}=e^1=e$. For several reformulations and applications of Schanuel's conjecture, see \cite{die2} and \cite{MS}.

Now, we shall find unconditional examples as corollaries of some known results from transcendental number theory. Let us mention them for making this text self-contained.

\begin{lemma}[Lindemann-Weierstrass]
Let $\,\alpha_1,\ldots,\alpha_n$ be distinct algebraic numbers. Then $e^{\alpha_1},\ldots,e^{\alpha_n}$ are linearly independent over $\QQ$.
\end{lemma}

For a proof of this theorem, see \cite[Theorem 1.4]{Baker}. 
\begin{corollary}
Let $\alpha_0,...,\alpha_n$ be nonzero algebraic numbers. If $f(z)=\sum_{i=0}^n\alpha_ie^{z+i}$, then $S_f=\emptyset$.
\end{corollary}
\begin{proof}
Let us take $z \in \QQ ~ \backslash \{0,-1,...,-n\}$, so that $0,z,z+1,...,z+k$ are distinct algebraic numbers. By Lindemann-Weierstrass theorem, the numbers $e^0,e^z,...,e^{z+n}$ are linearly independent over $\QQ$. Hence $f(z)$ is transcendental. In fact, towards a contradiction, suppose that $f(z)=\alpha \in \QQ$, the we would have the $\QQ$-linear relation
$$
\alpha e^0-\sum_{i=0}^n\alpha_ie^{z+i}=0.
$$
So, we need only to consider $z\in \{0,-1,...,-n\}$. In this case, $z+i=0$ for some $0\leq i\leq n$ and therefore if $f(z)=\alpha \in \QQ$, we get the $\QQ$-relation
$$
(\alpha-\alpha_z)e^0-\displaystyle\sum_{i=0}^n\alpha_ie^{z+i}=0
$$
This contradicts the fact that $e^0,e^{z}$

Hence $e^z + e^{z+1}$ has to be transcendental. For the omitted cases (i.e., $z=0$ and $z=-1$), $f(z)$ is transcendental because both $1+e$ and $e^{-1}+1$ are transcendental numbers. Therefore, $S_f = \emptyset$.
\end{proof}

Another important result comes from the work of Baker on linear forms of logarithms of algebraic numbers. It states that~\cite[Chap.~1]{Baker}:

\begin{lemma}[Baker]
Let $\,\alpha_1,\ldots,\alpha_n$ be non-zero algebraic numbers and let ~$\,\beta_0, \ldots, \beta_n$ be algebraic numbers, with $\beta_0 \neq 0$. Then the number $\,e^{\beta_0} \alpha_1^{\beta_1} \cdots \alpha_n^{\beta_n}$ is transcendental.
\end{lemma}

\begin{corollary}
If $g(z) = e^{\pi\,z+\alpha}$, where $\alpha\in \QQ\backslash \{0\}$, then $S_g=\emptyset$.
\end{corollary}
\begin{proof}
Suppose that $g(z)$ is algebraic for some $z \in \QQ$. Then by substituting $n=2$, $\alpha_1 = g(z) = e^{\pi\,z+\alpha}$, $\alpha_2 = -1$, $\beta_0 = -\alpha$, $\beta_1 = 1$, and $\beta_2 = i\,z$ in the Baker's theorem, with $i = \sqrt{-1}$, one finds that the number $e^{\beta_0} \, \alpha_1^{\beta_1} \, \alpha_2^{\beta_2}$ has to be transcendental. However, from the fact that $e^{i\,\pi} = -1$, one has
\begin{equation*}
e^{\beta_0} \, \alpha_1^{\beta_1} \, \alpha_2^{\beta_2} = e^{-\alpha} \left(e^{\alpha+ \pi\,z}\right)^1 \times (-1)^{i\,z} = e^{\pi\,z} \times \left(e^{i\,\pi}\right)^{i\,z} = 1 \, ,
\end{equation*}
which is an algebraic number, and we have a contradiction. Therefore, $g(z)$ can not be algebraic for any algebraic $z$ and then $S_g = \emptyset$.
\end{proof}


As our last tool, let us present the \emph{Mahler's classification scheme}: Let $P\in \Z[x]$ be a polynomial, the {\it height} of $P$, denoted by $\mathcal{H}(P)$, is the maximum of the absolute values of its coefficients. Let $\xi$ be a complex number, and for each pair of positive integers $n,h$, let $P(x)$ be a polynomial with degree at most $n$ and height at most $h$ for which $|P(\xi)|$ takes the smallest positive value and define $\omega(n,h)$ by the equation $|P(\xi)|=h^{-\omega(n,h)}$. Further, define
\begin{center}
$\omega_n=\limsup_{h\to \infty}\omega(n,h)$ and $\omega=\limsup_{n\to \infty}\omega_n$.
\end{center}

By defining $\nu(\xi)$ as the least positive integer $n$ for which $\omega_n(\xi)$ is infinity, we have four classes corresponding to the four possibilities for the values of $\omega(\xi)$ and $\nu(\xi)$:

\begin{itemize}
\item If $\omega(\xi)=0$, then $\xi$ is called an \textbf{$A$-number}.
\item If $0 < \omega(\xi) < \infty$, then $\xi$ is called an \textbf{$S$-number}.
\item If $\omega(\xi)=\infty$ and $\nu(\xi) < \infty$, then $\xi$ is called a \textbf{$U$-number}.
\item If $\omega(\xi)=\infty$ and $\nu(\xi) = \infty$, then $\xi$ is called a \textbf{$T$-number}.
\end{itemize}
Now, we state three useful properties of this classification:
\begin{itemize}
\item The set of $A$-numbers is precisely $\QQ$.
\item If $\alpha$ is a nonzero algebraic number, then $e^{\alpha}$ is an $S$-number.
\item If $\alpha$ and $\beta$ are two complex numbers having different Mahler classifications, then they are algebraically independent.
\end{itemize}

The proofs of these, and more, facts involving the Mahler's classification can be found in~\cite[Chap.~3]{Bugeaud}.

As an immediate consequence of the previous facts, we have
\begin{corollary}
Let $P(x,y)$ be a non-constant polynomial with algebraic coefficients and let $t$ be a $U$ or $T$-number. Then the function $h(z):=P(e^z,t)$ has an empty exceptional set.
\end{corollary}
\begin{proof}
Note that $h(0)$ is transcendental, because $h(0)=P(1,t)$ is transcendental (here we used the facts that $t$ is transcendental and $\QQ$ is algebraically closed). For $z\in \QQ\backslash \{0\}$, the number $e^z$ is a $S$-number and then if $P(e^z,t)=\alpha\in \QQ$, we would have an absurdity as $e^z$ and $t$ (a $U$ or $T$-number) algebraically dependent (since $Q(e^z,t)=0$, for $Q(x,y)=P(x,y)-\alpha$). Therefore $h(z)$ is transcendental. This completes the proof.
\end{proof}

\section{Proof of Theorem \ref{t2}}

In what follows, the symbol $x_1 \cdots \widehat{x_k} \cdots x_n$ will denote that $x_k$ is being omitted in the product $x_1 \cdots x_n$.

\begin{lemma}\label{1}
Given $n \geq 1$ complex numbers $a_0, \ldots, a_n$, not all zero, the polynomial
\begin{equation}
P_n(z) = a_0\,(z+1) \cdots (z^n + 1) + \sum_{k=1}^n{a_k \, z^k \, (z + 1) \cdots \widehat{(z^k + 1)} \cdots (z^n + 1)} \, ,
\end{equation}
is not the null polynomial.
\end{lemma}

\noindent
\textbf{Proof}
Suppose the contrary, i.e. $P_n(z)=0$, for all $z\in \C$. In particular, $P_n(0) = 0$ and then $a_0 = 0$. Thus
\begin{equation}
P_n(z) = \sum_{k=1}^n{a_k \, z^k \, (z + 1) \cdots \widehat{(z^k + 1)} \cdots (z^n + 1)} \, .
\end{equation}
Of course, its derivative $P_n'(z)$ is also null at $z=0$. On the other hand $P_n'(0) = a_1$. Therefore
\begin{equation}
P_n(z) = \sum_{k=2}^n{a_k \, z^k \, (z + 1) \cdots \widehat{(z^k + 1)} \cdots (z^n + 1)} \, .
\end{equation}
Similarly, we get $a_2 = a_3 = \cdots = a_n = 0$ from the nullity of higher derivatives $P_n''(0)$, $P_n^{(3)}(0), \ldots$, $P_n^{(n)}(0)$, respectively. But this contradicts our hypothesis on the coefficients $a_k$'s.
\qed
\bigskip

\begin{lemma}\label{2}
Let $P(z) \in \C[z]$ be a polynomial with degree $n$. Then
\begin{equation}
\label{desi1}
|P(z)| \leq \mathcal{L}(P) \, \max\{1,|z|\}^n \, ,
\end{equation}
for all $z \in \C$. Here $\mathcal{L}(P)$ is the \emph{length} of $P$, which is the sum of absolute values of its coefficients.
\end{lemma}

\noindent
\textbf{Proof}
If $P(z) = a_0 +a_1 z +\ldots +a_n z^n$, with $a_n \neq 0$, then
\begin{equation*}
|P(z)| = |a_0 +a_1 z +\ldots +a_n z^n| \leq |a_0| +|a_1|\,|z| +\ldots +|a_n|\,|z|^n \leq \mathcal{L}(P) \max\{1,|z|\}^n.
\end{equation*}
\qed
\bigskip

Recall that for $k\in \{1,...,n\}$, the function
\begin{equation*}
\sigma_k\left(x_1, \ldots, x_n \right):=  \displaystyle\sum_{1 \leq i_1< i_2 < \cdots < i_k \leq n} x_{i_1}x_{i_2} \cdots x_{i_k} 
\end{equation*}
is known as the \textit{elementary symmetric polynomial} in $x_1, \ldots, x_n$. An important property of these polynomials is that they satisfy the following identity:
\begin{equation}
\label{sim}
(z - x_1) \cdots (z - x_n) = z^n - \sigma_1 \left(x_1, \ldots, x_n \right)\,z^{n-1} +\ldots +(-1)^n \, \sigma_n \left(x_1, \ldots, x_n \right) .
\end{equation}

Now we can prove our main result.\\

\noindent
\textbf{Proof of Theorem \ref{t2}}
First, we claim that the function $\u$ is analytic in $\C^2$. In fact, set $Q_n(z) = \prod_{k=0}^n{(z-\alpha_k)}$. Then, by Lemma~\ref{2}, $|Q_n(z)| \leq \mathcal{L}(Q_n) \, \max\{1,|z|\}^n$, for all $z \in \C$. On the other hand, due to the identity in \eqref{sim}, $\mathcal{L}(Q_n) = 1 + \sum_{k=0}^n|\sigma_k(\alpha_1,\ldots,\alpha_n)|$. Thus
\begin{eqnarray*}
|\u(w,z)| & \leq & \sum_{n=1}^{\infty}\dfrac{|w|^n}{(1+\sum_{k=1}^n|\sigma_k(\alpha_1,\ldots,\alpha_n)|)(|w|^n+1)n!}\mathcal{L}(Q_n)\max\{1,|z|\}^n\\
 & \leq &  \sum_{n=1}^{\infty}\dfrac{\max\{1,|z|\}^n}{n!}=e^{\max\{1,|z|\}}
\end{eqnarray*}
Therefore, this function is well defined and, in fact, analytic in $\C^2$ (Weierstrass $M$-test). Given a pair $(t,\alpha)\in \R^{+} \times \QQ$, we have $\alpha = \alpha_{j+1}$ for some $j \geq 0$. Thus substituting $z = \alpha$ and $w = t$ in Eq.~\eqref{d} and multiplying both sides by $(t+1)(t^2+1) \cdots (t^j+1)$, we get
\begin{equation}
\label{desi2}
\sum_{k=1}^j{a_k t^k (t+1) \cdots \widehat{(t^k+1)} \cdots (t^j+1)} = \u(t,\alpha) \, \prod_{k=1}^j{\left(t^k+1 \right)}
\end{equation}
where the coefficients
\begin{equation*}
a_k = \frac{(\alpha-\alpha_1) \cdots (\alpha-\alpha_k)}{\left[1 +\sum_{k=1}^n{|\sigma_k(\alpha_1,\ldots,\alpha_n)|} \right] \, n!}
\end{equation*}
are algebraic numbers. According to Lemma~\ref{1}, Equation~\eqref{desi2} means that the number $t$ is a root of the non-zero polynomial
\begin{equation}
P_n(z) = -\u(t,\alpha)(z+1) \cdots (z^j + 1) +\sum_{k=1}^{j}{a_k z^k (z + 1) \cdots \widehat{(z^k + 1)} \cdots (z^j + 1)} \, .
\end{equation}
However the polynomial $P_n(z)\in \QQ[z]$ if and only if $\,\u(t,\alpha)\in \QQ$. Since the set $\QQ$ is algebraically closed and $P_n(t)=0$, we conclude that $t$ is transcendental if and only if $\,\u(t,\alpha)$ is transcendental.
\qed
\bigskip

%
%
%
%
%
%
%
%
%
%
%
%
%

%
%
\section*{Acknowledgements}

The authors would like to thank the anonymous referee for providing some important and enriching comments. The second author is also thankful to CNPq and FEMAT for financial support.
\footnotesize{

}

\end{document}